\documentclass[9pt,oneside,a4paper]{article}
\usepackage[margin=1.1in]{geometry}
\usepackage{amsmath,amsxtra,amssymb,latexsym, amscd,amsthm}
\usepackage{makeidx}
\usepackage{longtable}
\usepackage{multicol}
\usepackage{fancyhdr}
\usepackage{helvet} 
\usepackage{mathrsfs}
\usepackage[T1]{fontenc}
\usepackage{mwe}
\usepackage{subfig}
\usepackage{hyperref}
\usepackage{graphicx}
\usepackage{xcolor}
\usepackage{slashbox}
\usepackage{tikz}
\usepackage{pgfplots}

\theoremstyle{definition}
\newtheorem{define}{Definition}
\newtheorem{example}{Example}[section]
\theoremstyle{Plain}
\newtheorem{theorem}{Theorem}

\newtheorem{algorithm}{Algorithm}

\newtheorem{remark}{Remark}
\begin{document}
	
\title{Optimal Path Homotopy For Univariate Polynomials}


\author{Bao Duy Tran \\
	tranbaoduy@qnu.edu.vn\\
	Quy Nhon University, Vietnam}
\maketitle

\textbf{Abstract}: The goal of this paper is to study the path-following method for univariate polynomials. We propose to study the complexity and condition properties when the Newton method is applied as a correction operator. Then we study the geodesics and properties of the condition metric along those curves. Last, we compute approximations of geodesics and study how the condition number varies with the quality of the approximation. 

\section{Introduction}

The path-following method is known under several different names: homotopy method, numerical path-following, prediction-correction method, continuation method,... The history of homotopy methods is lengthier, which we will not try to discuss here. The basic idea is to find a path in the space of problems joining the problem we want to solve with a problem that is easy to solve and has the same structure. Then we follow the path starting from the easy problem and we use a correction operator to compute the solutions of the problems along the path. Finally, we end with an approximation of the solution to the problem we want to solve. This is a widely used method in optimization and system solving, as well as for interior-point methods. We refer to \cite{A94} for a general presentation. Here the problem we consider is root-finding for univariate polynomials. The complexity of this "king of methods" is linked to the conditioning of the polynomials along the chosen path, as shown in the work of Shub, Smale and many others. We refer to the survey by Dedieu \cite{D10} for this part. Recently, Beltran, Dedieu, Malajovich, and Shub studied geodesics for the condition metric and proved some convexity results for linear systems, see \cite{BDMS12}, log-convexity, and self-convexity on a Riemannian manifold, see \cite{BDMS10}. The interest of this condition metric is that the associated geodesics avoid ill-conditioned problems. In \cite{SS93, SS96}, the study of linear homotopy methods was in-depth; while later, in terms of the length of the path in the condition metric, a new bound was obtained for the complexity of path-following, see \cite{S09}.

This work focuses on some aspects of the path-following method applied to the particular case of univariate polynomial root-finding. In Section 2, we study geodesics from the definition and some examples and give the definitions of condition metrics. We provide some examples of the condition length of different curves in the space of polynomials and we derive conjectures on the approximation of geodesic with respect to the condition length. In this work, we want to approximate condition geodesics by B\'{e}zier curves. We, therefore, need to study and compute the condition length of a b\'{e}zier curve. To do the numerical computations in Matlab for the condition length of a B\'{e}zier curve, we define a general B\'{e}zier curve and its derivative. Next, in Section 3, we study the approximations of geodesics. We compute the length of condition geodesics approximated by B\'{e}zier curves in the space of univariate polynomials of finite degree. This is done through a minimization process. We study two cases: the space of univariate polynomials of degree $2$ and degree $3$. Finally, we consider the link between the complexity by the number of steps, required by the prediction-correction method, and the condition length of a curve to explain why it is interesting to study the condition metric.

\section{Geodesics and properties of the condition metric}
A geodesic is a locally length-minimizing curve. Geodesics depend on the chosen metric. 
\subsection{Geodesics and condition number}
There are many definitions of geodesics. We refer to \cite{C76} for the following presentation.
\begin{define}\label{def9}\cite[\text{Definition 8}, \text{p. 245}]{C76}
	A nonconstant, parameterized curve $\gamma:I\rightarrow S$ is said to be a geodesic at $t\in I$ if the field of its tangent vectors $\gamma'(t)$ is parallel along $\gamma$ at $t$, that is
	\[ \frac{D\gamma'(t)}{dt}=0, \]
	$\gamma$ is a parameterized geodesic if it is a geodesic for all $t\in I$.
\end{define}

The notion of a geodesic is local. Now we consider the definition of geodesic to subsets of $S$ that are regular curves.
\begin{define}\label{def10}\cite[\text{Definition 8a}, \text{p. 246}]{C76}
	A regular connected curve $C$ in $S$ is said to be a geodesic if, for every $p\in S$, the parameterization $\alpha(s)$ of a coordinate neighborhood of $p$ by the arc length $s$ is a parameterized geodesic; that is, $\alpha'(s)$ is parallel vector field along $\alpha(s)$.
\end{define}
\begin{example}\label{exa:3.1}\textit{Geodesics in Euclidean space}\\
	For the Euclidean metric, geodesics are line segments.
\end{example}

\begin{figure}[ht]
	\begin{minipage}[b]{0.6\linewidth}
		\begin{example}\label{exa:3.2}\textit{Geodesics on the sphere}\\
			The great circles of a sphere $S^2$ are geodesics. Indeed, the great circles $C$ are obtained by intersecting the sphere with a plane that passes through the center $O$ of the sphere. The principal normal at the point $p\in C$ lies in the direction of the line that connects $p$ to $O$ because $C$ is a circle of center $O$. Since $S^2$ is a sphere, the normal lies in the same direction, which verifies our assertion.
			\newline
			Through each point $p\in S$ and tangent to each direction in $T_p(S)$ that passes exactly one great circle, which is a geodesic. Therefore, by uniqueness, the great circles are the only geodesics of a sphere. 
		\end{example}	
	\end{minipage}
	\quad
	\begin{minipage}[b]{0.4\linewidth}
		\includegraphics[width=0.9\textwidth]{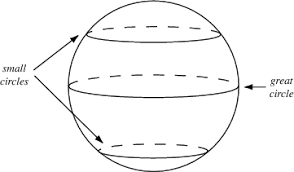}
		\caption{\small Geodesics on the sphere. Further information can be found on \cite{gc-w}.}\label{fig:1}
	\end{minipage}
	\normalsize
\end{figure}

\begin{figure}[ht]
	\begin{minipage}[b]{0.35\linewidth}
		\begin{example}\label{exa:3}\textit{Geodesics on the Poincar\'{e} plane}\\
			Consider the Poincar\'{e} plane. The sequence of arrows indicates how a tangent vector is rotated upon parallel transport along the curve. Vertical lines are geodesics, as are all semicircles that intersect the horizontal axis at a right angle.
		\end{example}	
	\end{minipage}
	\quad
	\begin{minipage}[b]{0.6\linewidth}
		\includegraphics[height=4cm,width=0.9\textwidth]{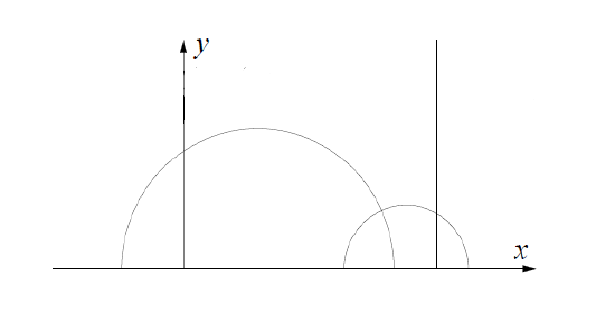}
		\caption{\small Geodesics on the Poincar\'{e} plane. Source: Author.}\label{fig:2}
	\end{minipage}
	\normalsize
\end{figure}
\newpage
To give the definitions of the condition metric specific to some different cases, we need to review some connotations of polynomial root-finding: condition number and iterative method.
\begin{define}\cite[\textit{Condition number}]{H02}\\
	Let $p(x)$ be a univariate polynomial of degree $d$ in the space of real polynomials: $p(x)=\sum\limits_{i=0}^{d}a_ix^i$. Suppose that $\alpha$ is a root of $p(x)$. Set $\tilde{p}(x)=\sum\limits_{i=0}^{d}\left| a_i\right| \left| x^i\right|$. We define the condition number for polynomial $p(x)$ as follows:
	\begin{equation}\label{eq:7}
	\mu(p,\alpha) = \frac{\left| \tilde{p}(\alpha)\right| }{\left| \alpha\right| \left| p^{'}(\alpha)\right|}.
	\end{equation}
\end{define}
\begin{remark}\label{rem1}\rm 
	To understand more about the condition number, we consider the simple case where $p(x) = x^2+bx+c$. This polynomial has two roots, denoted as $\alpha$ and $\beta$. Then we compute the condition number $(\ref{eq:7})$ and we obtain:
	\begin{equation*}
	\mu(p,\alpha) = \frac{\left| \alpha\right| ^2+\left| b\right| \left| \alpha\right| +\left| c\right| }{\left| \alpha\right| \left| 2\alpha+b\right|} \\
	= \frac{\left| \alpha\right| ^2+\left| -\alpha-\beta\right| \left| \alpha\right| +\left| \alpha\beta\right| }{\left| \alpha\right| \left| 2\alpha-\alpha-\beta\right| } \\
	= \frac{\left| \alpha\right| +\left| \alpha+\beta\right| +\left| \beta\right| }{\left| \alpha-\beta\right| }.
	\end{equation*}
	Since $\mu\rightarrow \infty$ when $\alpha=\beta$, we conclude that polynomials with close roots are ill-conditioned. 
\end{remark}
\begin{define}\textit{Newton method for univariate polynomials}\\
	Given a polynomial $p(x)\in\mathbb{R}\left[ x\right] $ and its derivative $p'(x)$, we begin with a first guess $x_0$ which is closed enough to a root of the polynomial $p$. The next iterate $x_1$ is defined as
	$$x_1=x_0-\dfrac{p(x_0)}{p^{'}(x_0)}.$$
	The process is repeated as
	$$x_{n+1}=x_n-\dfrac{p(x_n)}{p^{'}(x_n)}.$$
	until sufficient accuracy is reached.
\end{define}

The general idea is that the condition metric should be how large near the singular locus of the problem. When defining the condition metric, we essentially divide the Euclidean metric by the distance to the nearest singular system. See some examples below for more details.
\begin{example}\label{exa:3.4}
	Denote as $X=\mathbb{R}^2\setminus \left\lbrace\binom{0}{0}\right\rbrace$, consider $x=\binom{x_1}{x_2}$, $y=\binom{y_1}{y_2}$ . We consider the distance $d_1:X\times X\rightarrow \left[ 0,\infty\right)$ by
	\begin{equation*}
	d_1(x,y) = (x_1-y_1)^2+(x_2-y_2)^2,
	\end{equation*}
	The singular locus here is $\bar{0}=\binom{0}{0}$. The condition number of $x\in X$ is exactly the distance between $x$ and $\bar{0}$, given by
	\[ cond(x) = d_1(x,\bar{0})=x_1^2+x_2^2, \]
	The Euclidean metric $d:X\times X\rightarrow \left[ 0,\infty\right)$ is given by
	\[ d(x,y)=\sqrt{(x_1-y_1)^2+(x_2-y_2)^2}. \]
	Then we define the condition metric $d_k:X\times X\rightarrow \left[ 0,\infty\right)$ by
	\begin{equation*}
	d_k(x,y)=\frac{d(x,y)}{cond(x-y)}=\frac{\sqrt{(x_1-y_1)^2+(x_2-y_2)^2}}{(x_1-y_1)^2+(x_2-y_2)^2}=\frac{1}{\sqrt{(x_1-y_1)^2+(x_2-y_2)^2}}.
	\end{equation*}
	Therefore the condition norm on $X$ is given by
	\begin{equation}\label{eq:8}
	\left\| x\right\|_k=\frac{d(x,\bar{0})}{cond(x)}=\frac{1}{\sqrt{x_1^2+x_2^2}},
	\end{equation}
	(See Example \ref{exa:3.7} for using this norm.)
\end{example}
\begin{example}\label{exa:3.5}
	Denote as $\mathbb{R}\left[ x\right]_2$ the space of real univariate polynomials of degree $2$, consider $p=x^2+p_1x+p_2$ and $q=x^2+q_1x+q_2$. In term of vectors, $p=\binom{p_1}{p_2}$ and $q=\binom{q_1}{q_2}$. We consider the distance $d:\mathbb{R}\left[ x\right]_2\times \mathbb{R}\left[ x\right]_2\rightarrow \left[ 0,\infty\right) $ by
	\begin{equation*}
	d(p,q)=\left|\sqrt{\left(p_1-q_1 \right)^2-4\left| p_2-q_2\right|} \right|,
	\end{equation*}
	Then the condition number of $p$ is exactly the distance from $p$ to the singularity $\bar{0}=\binom{0}{0}$ on $\mathbb{R}\left[ x\right]_2$, is given by
	\begin{equation*}
	cond_c(p)=d(p,\bar{0})  = \left|\sqrt{p_1^2-4p_2} \right|.
	\end{equation*}
	The Euclidean metric $d:\mathbb{R}\left[ x\right]_2\times\mathbb{R}\left[ x\right]_2\rightarrow \left[ 0,\infty\right)$ is given by
	\[ d(p,q)=\sqrt{\left( p_1-q_1\right)^2+\left( p_2-q_2\right)^2},\]
	Then we can define the condition metric on $\mathbb{R}\left[ x\right]_2$ by
	\begin{equation*}
	d_c(p,q)=\frac{d(p,q)}{cond_c(p-q)}=\frac{\sqrt{\left( p_1-q_1\right)^2+\left( p_2-q_2\right)^2}}{\left|\sqrt{\left(p_1-q_1 \right)^2-4\left| p_2-q_2\right|} \right|}.
	\end{equation*}
	Hence, the condition norm on $\mathbb{R}\left[ x\right]_2$ is given by
	\begin{equation}\label{eq:9}
	\left\| p\right\|_c = \frac{d(p,\bar{0})}{cond_c(p)}=\frac{\sqrt{p_1^2+p_2^2}}{\left|\sqrt{p_1^2-4p_2} \right|}.
	\end{equation}
	(See Example \ref{exa:3.8} for using this norm.)
\end{example}
\begin{example}\label{exa:3.6}
	Denote as $\mathbb{C}\left[z\right]_n$ the space of complex univariate polynomials of degree $n$, consider $p(z)=z^{n}+a_{n-1}z^{n-1}+\ldots+a_1z+a_0$ and $q=z^n+b_{n-1}z^{n-1}+\ldots+b_1z+b_0$, where $a_i,b_i\in \mathbb{C}$, $i=0,1,\ldots,n-1$. In term of vectors,
	\[ p=\left[\begin{matrix}
	a_{n-1}, &
	a_{n-2}, &
	\ldots. &
	a_1, &
	a_0
	\end{matrix} \right]^T, \quad \text{and} \quad q=\left[\begin{matrix}
	b_{n-1}, &
	b_{n-2}, &
	\ldots, &
	b_1, &
	b_0
	\end{matrix} \right]^T. \]
	The Euclidean metric $d:\mathbb{C}\left[ z\right]_n\times\mathbb{C}\left[ z\right]_n\rightarrow \left[ 0,\infty\right)$ is given by
	\[ d(p,q)=\sqrt{\sum\limits_{i=0}^{n-1}\left( a_i-b_i\right)^2}.\]
	Then we can define the condition metric on $\mathbb{C}\left[z\right]_n$ by
	\[ d_{cn}(p,q)=\frac{d(p,q)}{cond_{cn}(p-q)}=\frac{\sqrt{\sum\limits_{i=0}^{n-1}\left( a_i-b_i\right)^2}}{\left|\sqrt[n]{D(p-q)}\right|}, \]
	where $D(p-q)$ is the discriminant of $p-q$, see its definition in \cite[Chapter 1]{B08}.\\
	The condition norm on $\mathbb{C}\left[z\right]_n$ is given by
	\begin{equation}\label{eq:10}
	\left\| p \right\|_{cn}=\frac{d(p,\bar{0})}{cond_{cn}(p)}=\frac{\left(\sum\limits_{i=0}^{n-1}a_i^2\right)^{\frac{1}{2}}}{\left|\sqrt[n]{D(p)}\right| }.
	\end{equation}
	where $\bar{0}=\left[ 0\; 0\; \ldots \;0\; 0\right]^T$ and $D(p)$ is the discriminant of $p$.\\
	(See Example \ref{exa:3.9} for using this norm.)
\end{example}
\subsection{Some examples of condition length and conjecture geodesics}
Now we study some examples to observe how to conjecture geodesics. In which, the condition length holds an important role.
\begin{example}\label{exa:3.7}
	We give a toy example for the geodesics of condition metric. The space of systems is identified with $\mathbb{R}^2$ (as in Example \ref{exa:3.4}) and the singular systems are reduced to the point $\binom{0}{0}$ and hence the condition number of the system $\binom{x}{y}$ is $\frac{1}{x^2+y^2}$. So what are the geodesics in this case?\\
	Take $A,B\in\mathbb{R}^2\setminus \left\lbrace \binom{0}{0}\right\rbrace$, let $\gamma\in \mathcal{C}\left( \left[ 0,1\right],\mathbb{R}^2\setminus\left\lbrace \binom{0}{0} \right\rbrace \right)$ be a path between $A$ and $B$ avoiding $\binom{0}{0}$. The length of $\gamma$ for the metric is
	\[ l_k(\gamma)=\int\limits_{0}^{1}\frac{dt}{\left\| \gamma(t)\right\|_k^2}. \]
	Consider a simple case where $A=\binom{-1}{0}$, $B=\binom{1}{0}$ and $C=\binom{0}{1}$ (see Fig. \ref{fig:3}). First, we want to compute the condition length of $\gamma$.
	\newline
	We consider two paths. One is the piecewise segment path
	$$ \gamma_1(t) = \begin{cases}
	\gamma_{1,1}(t) = \binom{2t-1}{2t}, & t\in \left[ 0,\frac{1}{2}\right],\\
	 \gamma_{2,1}(t) = \binom{2t-1}{2t-2}, & t\in \left[ \frac{1}{2},1\right].
	\end{cases} $$
	The second one is the arc of unit circle going through $C$
	\[ \gamma_2(t)= e^{i\pi(t-1)}=\binom{\cos(\pi(t-1))}{\sin(\pi(t-1))},\quad t\in \left[0,1\right]. \]
	Then we have,
	\[l_c(\gamma_1)=\int\limits_{0}^{1}\dfrac{dt}{\left\| \gamma_1(t)\right\|_k^2}=\dfrac{\pi}{2}, \text{ and } 	l_c(\gamma_2) = 1. \]

	So for the condition metric, we conclude that 
	\[ l_c(\gamma_1)=\frac{\pi}{2}\simeq 1.6 \geq 1 = l_c(\gamma_2). \]
	But for the Euclidean metric, it's easy to see that
	$l(\gamma_1)=2\sqrt{2}\simeq 2.83 \leq 3.14\simeq\pi=l(\gamma_2)$. This is one of the reasons why we are interested in the condition metric.
	
	\begin{figure}[ht]
		\centering
	\begin{tikzpicture}
	\draw[->] (-2.5,0) -- (3,0) node[right] {$x$};
	\draw[-] (0,-2) -- (0,0);
	\draw[dashed,->] (0,0) -- (0,4.2) node[above] {$\frac{1}{x^2+y^2}$};
	\draw[-] (-1.5,-1.5) -- (0,0);
	\draw[dashed,-] (0,0) -- (1,1);
	\draw[->] (1,1) -- (2,2) node[right] {$y$};
	\draw[scale=1,domain=0.5:2,smooth,variable=\x,blue] plot ({\x},{1/(\x*\x)});
	\draw[scale=1,domain=0.5:2,smooth,variable=\x,blue] plot ({-\x},{1/(\x*\x)}) node[above] {(S)};
	
	\draw[scale=1,domain=-1:1,smooth,variable=\x,red]  plot ({\x},{-sqrt(1-\x*\x)});
	\draw[scale=1,domain=-1:0,smooth,variable=\x,green] plot ({\x},{-\x - 1});
	\draw[scale=1,domain=0:1,smooth,variable=\x,green] plot ({\x},{\x - 1});
	
	\filldraw (-1,0) circle[radius=0.5pt] node[below left] {A};
	\filldraw (1,0) circle[radius=0.5pt] node[below right] {B};
	\filldraw (0,-1) circle[radius=0.5pt] node[below right] {C};
	
	\draw[scale=1,domain=-0.707:0.707,smooth,variable=\x,red]  plot ({\x},{2-sqrt(0.5-\x*\x)});
	\draw[scale=1,smooth, variable=\x,green] (-0.0175,1.285) arc (0:90:0.7);
	\draw[scale=1,smooth, variable=\x,green] (0.725,2) arc (90:180:0.7);
	
	\draw[scale=1,domain=-0.55:0.55,smooth,variable=\x, blue]  plot ({\x},{3.25-sqrt(0.325-\x*\x)});
	\draw[dashed,scale=1,domain=-0.55:0.55,smooth,variable=\x, blue]  plot ({\x},{3.875-sqrt(0.325+\x*\x)});
	\draw[scale=1,domain=-1.175:1.175,smooth,variable=\x, blue]  plot ({\x},{1.55-sqrt(2-\x*\x)});
	\draw[dashed,scale=1,domain=-1.175:1.175,smooth,variable=\x, blue]  plot ({\x},{2.575-sqrt(2+\x*\x)});

	\draw[green,-] (2.75,4) -- (3,4) node[right] {$\gamma_1$};	
	\draw[red,-] (2.75,3.5) -- (3,3.5) node[right] {$\gamma_2$};
	\end{tikzpicture}
	\caption{\small Geodesic of condition metric and its "lifting" visualization. Source: Author.}
	\label{fig:3}
	\normalsize
	\end{figure}
	
	In other to have a better visualize of this example, we consider the surface $(S)$ parameterized by $\left(x,y,\frac{1}{x^2+y^2} \right)$ in $\mathbb{R}^3$. Then we can "lift" the paths $\gamma_1$ and $\gamma_2$ on $(S)$ as on Fig. \ref{fig:3}, which were plotted by TikZ package in LaTeX. In fact, "lifting on $(S)$" sends the condition metric in $\mathbb{R}^2$ to the metric induced on $(S)$ by the Euclidean metric in $\mathbb{R}^3$.\\
	With this interpretation, we can see the clear difference in length between the two paths. Therefore we conjecture that $\gamma_2$ is a geodesic for the condition metric.
\end{example}

\begin{example}\label{exa:3.8}
	In the space of real polynomials of degree $2$, let
	$p_1(x)=x^2+b_1x+c_1$, and 	$p_2(x)=x^2+b_2x+c_2$.
	Given a path joining two polynomials, we want to compute its condition length.
	\newline
	The linear homotopy path between $p_1$ and $p_2$ is given by
	\begin{equation}\label{eq:11}
	\gamma(t,x)=(1-t)p_1(x)+tp_2(x).
	\end{equation}
	We see that each polynomial $p(x)=x^2+bx+c$ can be seen as a point $(b,c)\in\mathbb{R}^2$.
	So we just identify $p_1(x)$ to $(b_1,c_1)$ and $p_2(x)$ to $(b_2,c_2)$.
	\newline
	The linear homotopy (\ref{eq:11}) can be written as $\gamma(t,x) = x^2+\left( (1-t)b_1+tb_2 \right)x + \left( (1-t)c_1+tc_2 \right)$\\
	or, as a vector in $\mathbb{R}^2$
	\[\gamma_t = \left( (1-t)b_1+tb_2,(1-t)c_1+tc_2 \right).\]
	(It is represented by the segment of line joining $p_1$ and $p_2$).\\
	The condition length of segment $\overline{p_1p_2}$ is equal to
	\begin{align*}
	l_c(\overline{p_1p_2}) &= \int\limits_{0}^{1} \left\| \dot{\gamma}(t)\right\|_c dt \\
	&= \int\limits_{0}^{1} \left\| \binom{b_2-b_1}{c_2-c_1} \right\| \dfrac{dt}{\left| \sqrt{\left( (1-t)b_1+tb_2\right)^2-4\left( (1-t)c_1+tc_2 \right)  }\right|  } \\
	&= \sqrt{(b_2-b_1)^2+(c_2-c_1)^2} \int\limits_{0}^{1}\dfrac{dt}{\left| \sqrt{\left( (1-t)b_1+tb_2\right)^2-4\left( (1-t)c_1+tc_2 \right)}\right|}.
	\end{align*}
		
		\begin{figure}[ht]
		\begin{minipage}[b]{0.55\linewidth}
			For instance, we want to compute the condition length path between $p(x)=x^2-x-1$ and $q(x)=x^2+x-1$.
			In term of vectors (in $\mathbb{R}^2$), $p=\binom{-1}{-1}$ and $q=\binom{1}{-1}$. Let us consider two paths. The first one is an arc of circle going through $r=\binom{0}{-2}$ that is parameterized by
			\[ \gamma_1(t)=\binom{\cos\left(\pi(t-1)\right) }{-1+\sin\left(\pi(t-1) \right) },\qquad t\in\left[ 0,1\right]. \]
			The second one is a piecewise segment path going through $p$, $r$ and $q$, which is given by
			\begin{align*}
			\gamma_2(t)&=\left\lbrace \begin{array}{l}
			\gamma_{1,2}(t) = (1-2t)p+2tr, \quad\quad t\in\left[ 0,\frac{1}{2}\right], \\
			\gamma_{2,2}(t) = (2-2t)r+(2t-1)q, \; t\in\left[\frac{1}{2},1\right], \\
			\end{array} \right. \\
			&= \left\lbrace \begin{array}{l}
			\gamma_{1,2}(t)=\binom{2t-1}{1-2t},\quad t\in\left[ 0,\frac{1}{2}\right], \\
			\gamma_{2,2}(t)=\binom{2t-1}{2t-3},\quad t\in\left[\frac{1}{2},1\right].\\
			\end{array} \right. 
			\end{align*}
		\end{minipage}
		\quad
		\begin{minipage}[b]{0.4\linewidth}
			\includegraphics[height=7cm,width=0.9\textwidth]{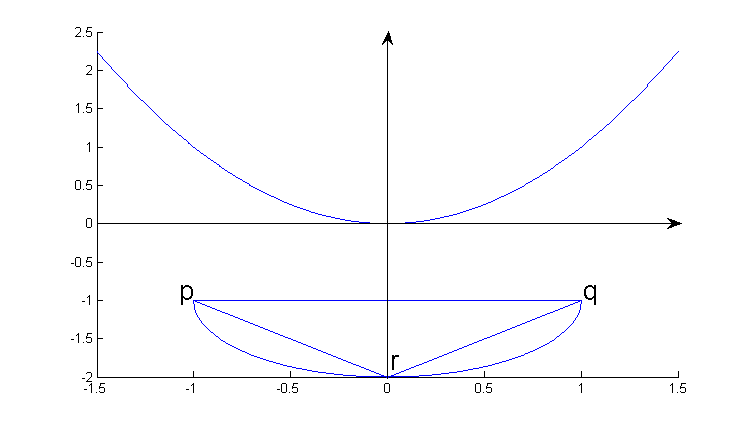}
			\caption{\small Example \ref{exa:3.8}. Source: Author.}
		\end{minipage}
		\normalsize
		\end{figure}	
	Thus, we can get the condition length for $\gamma_1$ and $\gamma_2$ as follows:
	\begin{align*}
	l_c(\gamma_1)&=\int\limits_{0}^{1}\left\| \dot{\gamma_1}(t)\right\|_c dt = \pi\int\limits_{0}^{1} \dfrac{dt}{\left|\sqrt{\cos^2\left(\pi(t-1) \right)-4\left( -1+\sin\left(\pi(t-1) \right)  \right)}\right| } \simeq 1.191\\
	l_c(\gamma_2) &= \int\limits_{0}^{\frac{1}{2}} \left\| \dot{\gamma}_{1,2}(t)\right\| _c dt + \int\limits_{\frac{1}{2}}^{1}\left\| \dot{\gamma}_{2,2}(t)\right\| _c dt \\
	&= 2\sqrt{2}\int\limits_{0}^{\frac{1}{2}}\dfrac{dt}{\left| \sqrt{(2t-1)^2-4(1-2t)}\right| }+2\sqrt{2}\int\limits_{\frac{1}{2}}^{1}\dfrac{dt}{\left| \sqrt{(2t-1)^2-4(2t-3)}\right| }\\
	&\simeq 1.586.
	\end{align*}
	So for the condition metric, we conclude that
	\[ l_c(\gamma_2) \geq l_c(\gamma_1), \] 
	By the same way and similar interpretation as at the end of Example \ref{exa:3.7}, we conjecture that $\gamma_1$ is a geodesic for the condition metric.
\end{example}

Now we take into account a more general case. 
\begin{example}\label{exa:3.9}
	In the space of complex polynomials of degree $n$ ($n\geq 2$), let
	$p(z)=z^n+a_{n-1}z^{n-1}+\ldots+a_1z+a_0$, and $q(z)=z^n+b_{n-1}z^{n-1}+\ldots+b_1z+b_0$.
	Given a path joining two polynomials, we want to compute its condition length.\\
	In term of vectors (in $\mathbb{C}^n$), $p=\left[a_{n-1}\; a_{n-2}\; \ldots \; a_1 \; a_0\right]^T $ and $q=\left[ b_{n-1}\; b_{n-2}\; \ldots \; b_1 \; b_0\right]^T$.\\
	The linear homotopy path between $p$ and $q$ is given by
	\begin{align*}
	f(t) &= (1-t)p+tq, \quad t\in\left[ 0,1\right] \\
	&= \left[\begin{matrix}
	(1-t)a_{n-1}+tb_{n-1}, &
	(1-t)a_{n-2}+tb_{n-2}, &
	\ldots &
	(1-t)a_1+tb_1, &
	(1-t)a_0+tb_0
	\end{matrix} \right]^T, \; t\in\left[ 0,1\right].
	\end{align*}
	Then we have
	\[ f'(t)=\left[\begin{matrix}
	b_{n-1}-a_{n-1}, &
	b_{n-2}-a_{n-2}, &
	\ldots &
	b_{1}-a_{1}, &
	b_{0}-a_{0}
	\end{matrix} \right]^T, \; t\in \left[ 0,1\right].  \]
	The condition length of $f$ is equal to
	\begin{equation*}
	I=\int\limits_{0}^{1} \frac{\left\| f'(t)\right\|}{cond_{cn}\left(f(t)\right)} dt= \int\limits_{0}^{1} \frac{\sqrt{\sum\limits_{i=0}^{n-1}\left(b_i-a_i\right)^2}}{\left|\sqrt[n]{D(f)} \right|}dt,
	\end{equation*}
	where $D(.)$ is the discriminant of $f$.
\end{example}

The approximation of geodesic needs a "material" called B\'{e}zier curves. We will give their definition and some properties in the following subsection.
\subsection{B\'{e}zier curves and its condition lengths}
A B\'{e}zier curve is defined by a set of \textit{control points} $P_0, P_1,\ldots, P_d$, where $d$ is called the degree of the curve ($d=1$ for linear, $d=2$ for quadratic, $d=3$ for cubic, etc.). Denotes $B\left( \left[ P_0,P_1,\ldots,P_d\right] ,t\right) $ the parameterization of B\'{e}zier curve of degree $d$, $d=0,1,\ldots, d$, associated with $P_0,P_1,\ldots,P_d$. B\'{e}zier curves can be defined for any degree $d$. A B\'{e}zier curve of degree $d$ is a point-to-point linear combination of a pair of corresponding points in two B\'{e}zier curves of degree $d-1$
	\begin{equation}\label{eq:12}
	B\left( \left[P_0,P_1,\ldots,P_d \right],t \right)
	=(1-t)B\left( \left[P_0,P_1,\ldots,P_{d-1} \right],t \right)+tB\left( \left[P_1,P_2,\ldots,P_d \right],t \right), \; t\in\left[ 0,1\right].
	\end{equation}
Then, since $b_{i,d}(t)$, where $d\in \mathbb{N}$ and $i\in \{0,\ldots,d \}$, are the Bernstein polynomials, $\{b_{0,d},b_{1,d}, \ldots, b_{d,d}\}$ is the basis of $\mathbb{R}[t]_d$ and every polynomial parametrization of a curve can be seen as a B\'{e}zier parametrization, see \cite[Chapter 1]{B08}, the formula (\ref{eq:12}) can be expressed explicitly as follows:
				\begin{equation}\label{eq:13}
				B\left( \left[P_0,P_1,\ldots,P_d \right],t \right)= \sum\limits_{i=0}^{d}P_ib_{i,d}(t),
				\end{equation}
				where $b_{i,d}(t)$, $i=0,1,\ldots,d$, are the Bernstein polynomials defined in \cite[Chapter 1]{B08}.
\begin{remark}\label{rem3}\rm 
	Based on the formula $(\ref{eq:13})$, the derivative for a B\'{e}zier curve of degree $d$ is given by
	\begin{equation}\label{eq:14}
	B'(t)=\frac{d}{dt}B([P_0,P_1,\ldots,P_d],t)=d\sum\limits_{i=0}^{d-1}b_{i,d-1}(t)\left(P_{i+1}-P_i \right).
	\end{equation}
\end{remark}
From formula (\ref{eq:13}) and Eq. (\ref{eq:14}), we could obtain the parameterization of a B\'{e}zier curve and its derivative. Therefore, we can compute its condition length (in the space of polynomials of degree $2$, using condition norm (\ref{eq:9})).\\
For a more general setting, by using condition norm (\ref{eq:10}), we can compute the condition length of a B\'{e}zier curve in the space of univariate polynomials of degree $n$, where $n\geq 2$.
\subsection{Properties of the condition length}
Let us now examine at the characteristics of the condition length. \\
Consider the B\'{e}zier curve
\[ \Gamma(t)=B\left(\left[ P_0,P_1,\ldots,P_d\right],t\right)=\sum\limits_{i=1}^d\binom{d}{i}(1-t)^{d-i}t^iP_i,\qquad t\in\left[ 0,1\right]. \]
where $P_i\in\mathbb{C}\left[ z\right]_d$, $i=0,1,\ldots,d$, i.e. each control point $P_i$, $i=0,1,\ldots,d$, is a degree $n$ complex polynomial. Suppose that $P_i=z^n+P_{i,n-1}z^{n-1}+\ldots+P_{i,1}z+P_{i,0}$, $i=0,1,\ldots,d$, i.e. in term of vectors, we consider
\[ P_i=\left[\begin{matrix}
P_{i,n-1} &
P_{i,n-2} &
\ldots &
P_{i,1} &
P_{i,0}
\end{matrix} \right]^T. \]
The condition length of $\Gamma$ is given by the function $\texttt{lc}\left(\Gamma\right)$ as below:
\begin{equation}\label{eq:15}
\texttt{lc}\left(\Gamma\right)=\texttt{lc}\left(P_0,P_1,\ldots,P_d\right)=\int\limits_{0}^{1}\left\| \Gamma'(t)\right\|_{cn}dt=\int\limits_{0}^{1}\frac{\left\| \Gamma'(t)\right\|}{cond_{cn}\left(\Gamma(t)\right)}dt.
\end{equation}
\begin{remark}\label{remA}\rm 
	We have
	\[ \left\| \Gamma'(t)\right\|=\left\| d\sum\limits_{i=0}^{d-1}\binom{d-1}{i}(1-t)^{d-1-i}t^i\left(P_{i+1}-P_i\right)\right\|.  \]
	We see that $\Gamma'(t)$ is a differentiable function with respect to the coordinates of $P_0,P_1,\ldots,P_d$, and $\Gamma'(t)\neq 0$, so the composed function $\left\| \Gamma'(t)\right\|$ is also differentiable.\\
	Furthermore,
	\[ cond_{cn}\left(\Gamma(t)\right)=cond_{cn}\left(\sum\limits_{i=0}^{d}\binom{d}{i}(1-t)^{d-i}t^iP_i\right),  \]
	is a non-zero differentiable function.\\
	As a result, $\texttt{lc}$ (as defined in ($\ref{eq:15}$)) is a differentiable function with respect to the coordinates of $P_0,P_1,\ldots,P_d$.
\end{remark}
\begin{remark}\label{remB}\rm 
	Moreover, for any $i=0,1,\ldots,d$; $j=0,1,\ldots,n-1$, we have
	\begin{equation*}
	\frac{\partial}{\partial P_{i,j}}\texttt{lc}\left(\Gamma\right)=\frac{\partial}{\partial P_{i,j}}\int\limits_{0}^{1}\frac{\left\| \Gamma'(t)\right\|}{cond_{cn}\left(\Gamma(t)\right)}dt=\int\limits_{0}^{1}\frac{\partial}{\partial P_{i,j}}\left( \frac{\left\| \Gamma'(t)\right\|}{cond_{cn}\left(\Gamma(t)\right)}\right)dt.
	\end{equation*}
	We see that $\frac{\partial}{\partial P_{i,j}}\left( \frac{\left\| \Gamma'(t)\right\|}{cond_{cn}\left(\Gamma(t)\right)}\right)$, for any $i=0,1,\ldots,d$, $j=0,1,\ldots,n-1$, is a rational function where numerator and denominator are differentiable functions, hence $\frac{\partial}{\partial P_{i,j}}\left( \frac{\left\| \Gamma'(t)\right\|}{cond_{cn}\left(\Gamma(t)\right)}\right)$ is a differentiable function for any $i=0,1,\ldots,d$, $j=0,1,\ldots,n-1$.\\
	Therefore $\frac{\partial}{\partial P_{i,j}}\texttt{lc}\left(\Gamma\right)$ is a differentiable function for any $i=0,1,\ldots,d$, $j=0,1,\ldots,n-1$, i.e, $\texttt{lc}$ (as defined in ($\ref{eq:15}$)) is a twice differentiable function with respect to coordinates of $P_0,P_1,\ldots,P_d$.
\end{remark}
Remarks \ref{remA} and \ref{remB} allow us to compute the first and second derivatives of the condition length of a B\'{e}zier curve in $\mathbb{C}\left[ z\right]_d$, where $d\geq 2$. We will use their derivatives in the next section. 

\section{Approximations of geodesics}
We notify that all computations are done by using Matlab, and tests were performed on a machine running \texttt{Windows 10 Pro} with an \texttt{Intel(R) Core(TM) i5-8265U \@ 1.6GHz 1.8GHz} and \texttt{Installed RAM 8GB}. \\
We will approximate the geodesics by the B\'{e}zier curve. Our study will place in the space of univariate polynomials of degree $n$, $n\in \mathbb{N}$, $n\geq 2$. We use the Euclidean norm $\left\| .\right\|$ with the appropriate dimension for the following calculation.

\subsection{In the space of univariate polynomials of degree 2}
To understand more about the condition length of the B\'{e}zier curve, we study some examples below, where we use the B\'{e}zier curve defined by these control points as an initial guess for an optimization process that approximates the condition geodesic of endpoints $\binom{-1}{-1}$ and $\binom{1}{-1}$. Note that the endpoints are kept fixed throughout the optimization process.\\
In this subsection, \texttt{lc(p)} indicates the condition length of the B\'{e}zier curve defined by \texttt{p}. We obtain the optimal value, together with the matrix of control points (\texttt{x\_opt}) for the B\'{e}zier curve that realizes the minimum and its condition length (\texttt{fval}).
\begin{example}\label{exa:4.1}
	Consider the set of control points $\left\lbrace \binom{-1}{-1},\binom{0}{-2},\binom{1}{-1}\right\rbrace$ which we denote using a matrix:
	\small
	$$\texttt{p} = \begin{pmatrix}
		-1.0000 & 0 & 1.0000\\
		-1.0000 & -2.0000 & -1.0000
	\end{pmatrix},$$
	\normalsize
	therein we mean that these control points for the control polygon of a B\'{e}zier curve of degree $2$. \\
	We get the condition length $\texttt{lc(p)} = 1.4524$, and the optimal value
	$$ \texttt{fval} = 1.3948, \quad \texttt{x\_opt} = \begin{pmatrix}
	 -1.0000 & 0 & 1.0000 \\
	-1.0000 & -1.4407 & -1.0000
	\end{pmatrix}.$$
	Now we are going to perturb the optimal control points a little to see how large is the variation of the optimal condition length. Let's consider other sets of control points $\left\lbrace \binom{-1}{-1},\binom{0}{-1.4},\binom{1}{-1}\right\rbrace$ and $\left\lbrace \binom{-1}{-1},\binom{0}{-1.5},\binom{1}{-1}\right\rbrace$, which we denote using matrices \texttt{p1} and \texttt{p2} as below:
	\small
	$$\texttt{p1} = \begin{pmatrix}
	-1.0000 & 0 & 1.0000\\
	-1.000 & -1.4000 & -1.0000
	\end{pmatrix}, \qquad \texttt{p2} = \begin{pmatrix}
	-1.0000 & 0 & 1.0000\\
	-1.0000 & -1.5000 & -1.0000
	\end{pmatrix}.$$
	\normalsize
\begin{table}[htb]
	\centering
	\begin{tabular}{c  c  c  c  c  c }
		\hline
		  \texttt{lc(p1)} & \texttt{lc(p2)} & $\left\| \texttt{p1}-\texttt{x\_opt} \right\| $ & $\left\| \texttt{p2}-\texttt{x\_opt} \right\| $ & $\left\| \texttt{lc(p1)}-\texttt{fval} \right\| $ & $\left\| \texttt{lc(p2)}-\texttt{fval} \right\| $  \\
		1.3952 & 1.3956 & 0.0407 & 0.0593 & 3.7227e-004 & 7.6385e-004 \\
		\hline
		\hline
		\texttt{lc(q1)} & \texttt{lc(q2)} & $\left\| \texttt{p}-\texttt{q1} \right\| $ & $\left\| \texttt{p}-\texttt{q2} \right\| $ & $\left\| \texttt{lc(p)}-\texttt{lc(q1)} \right\| $ & $\left\| \texttt{lc(p)}-\texttt{lc(q2)} \right\| $  \\
		1.4350 & 1.4721 & 0.1000 & 0.1000 & 0.0174 & 0.0197 \\
		\hline
	\end{tabular}
	\caption{\small The condition lengths of the B\'{e}zier curves are defined by \texttt{p1, p2}, \texttt{q1, q2}, and the perturbations with their optimum.}\label{tab3.1}
	\normalsize
\end{table}
\begin{figure}[!ht]
	\subfloat[Starting control points \texttt{p, q1, q2}, and its B\'{e}zier curve.\label{subfig-1:5}]{%
		\includegraphics[width=0.45\textwidth]{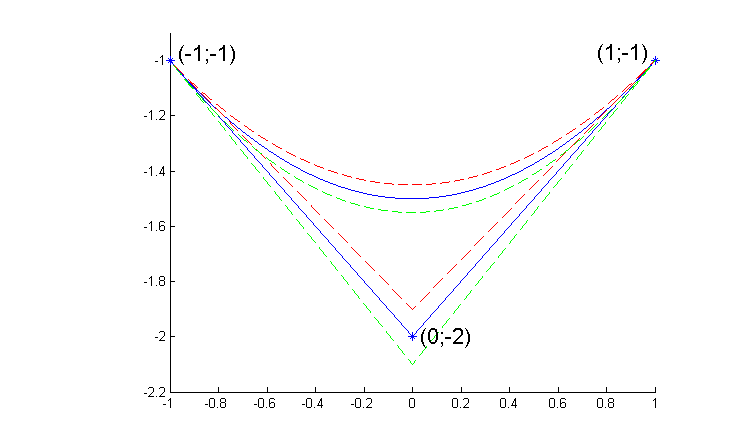}
	}
	\hfill
	\subfloat[The optimal control points \texttt{x\_opt} and its approximate geodesic.\label{subfig-2:5}]{%
		\includegraphics[width=0.45\textwidth]{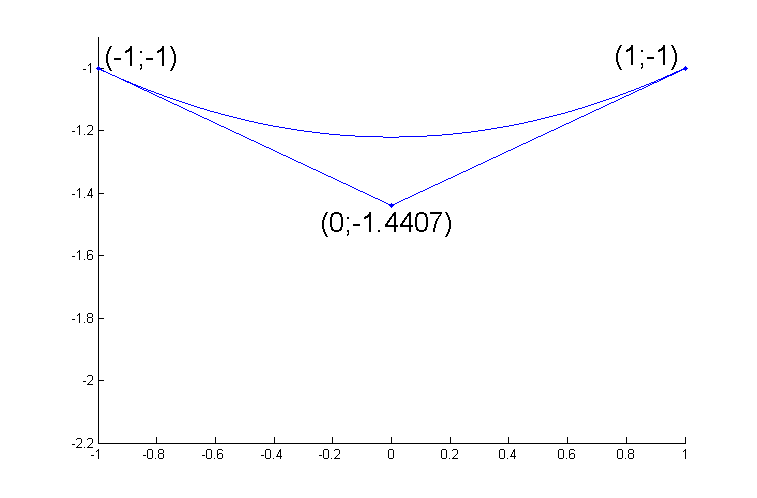}
	}
	\caption{\small B\'{e}zier curve of degree 2 that approximates the condition geodesic of endpoints $\binom{-1}{-1}$ and $\binom{1}{-1}$. Source: Author.}
	\label{fig:5}
	\normalsize
\end{figure}
	By the above estimation $\left\| \texttt{lc(p1)}-\texttt{fval} \right\| $ and $\left\| \texttt{lc(p2)}-\texttt{fval} \right\|$ in Table \ref{tab3.1}, we can conclude the optimal value of the condition length does not change much, with the gap values less than $0.001$.
	
	On the other hand, we will perturb the initial control points a bit to see how large the condition length variation is. We consider other sets of control points $\left\lbrace \binom{-1}{-1},\binom{0}{-1.9},\binom{1}{-1}\right\rbrace$ and $\left\lbrace \binom{-1}{-1},\binom{0}{-2.1},\binom{1}{-1}\right\rbrace$, which we denote using matrices \texttt{q1} and \texttt{q2} as follows:
	$$\texttt{q1} = \begin{pmatrix}
	-1.0000 & 0 & 1.0000\\
	-1.0000 & -1.9000 & -1.0000
	\end{pmatrix}, \qquad \texttt{q2} = \begin{pmatrix}
	-1.0000 & 0 & 1.0000\\
	-1.0000 & -2.1000 & -1.0000
	\end{pmatrix}.$$
	By taking the perturbation as in the two bottom lines of Table \ref{tab3.1}, we can conclude the condition length does not change much, with the gap values less than $0.0198$.
\end{example}
\begin{example}\label{exa:4.2}
	We are going to use a B\'{e}zier curve of degree $3$. Consider the set of control points $\left\lbrace \binom{-1}{-1}, \binom{-0.5}{-2}, \binom{0.5}{-2.5},\binom{1}{-1}\right\rbrace $ which we denote using a matrix:
	$$ \texttt{p} = \begin{pmatrix}
		-1.0000 & -0.5000 & 0.5000 & 1.0000 \\
		-1.0000 & -2.0000 & -2.5000 & -1.0000
	\end{pmatrix},$$
	here we mean that these control points for the control polygon of a B\'{e}zier curve of degree $3$.\\
	The condition length is computed $\texttt{lc(p)} = 1.5090$. Then we approximate a condition geodesic joining $\binom{-1}{-1}$ and $\binom{1}{-1}$, we obtain the optimal value of the geodesic as below
	$$ \texttt{fval} = 1.2398, \quad \texttt{x\_opt} = \begin{pmatrix}
	-1.0000 & -0.3753  &  0.3753 & 1.0000 \\
	-1.0000 & -1.2910 & -1.2909 & -1.0000
	\end{pmatrix}.$$
	Now we are going to perturb the optimal control points a little to see how large is the variation of the optimal condition length. Let's consider other sets of control points $\left\lbrace \binom{-1}{-1},\binom{-0.3753}{-1.2},\binom{0.3753}{-1.2},\binom{1}{-1}\right\rbrace $, $\left\lbrace \binom{-1}{-1}, \binom{-0.3753}{-1.35}, \binom{0.3753}{-1.35},\binom{1}{-1}\right\rbrace $ and $\left\lbrace \binom{-1}{-1},\binom{-0.4}{-1.25},\binom{0.38}{-1.26},\binom{1}{-1}\right\rbrace $, which we denote using matrices \texttt{p1}, \texttt{p2} and \texttt{p3} as below
	\small
	$$\texttt{p1} = \begin{pmatrix}
	-1.0000 & -0.3753  &  0.3753 & 1.0000 \\
	-1.0000 & -1.2000 & -1.2000 & -1.0000
	\end{pmatrix},\; \texttt{p2} = \begin{pmatrix}
	-1.0000 & -0.3753  &  0.3753 & 1.0000 \\
	-1.0000 & -1.3500 & -1.3500 & -1.0000
	\end{pmatrix},$$ 
	 $$\texttt{p3} = \begin{pmatrix}
	-1.0000 & -0.4000  &  0.3800 & 1.0000 \\
	-1.0000 & -1.2500 & -1.2600 & -1.0000
	\end{pmatrix}.$$
	\normalsize
	\begin{table}[htb]
		\footnotesize
		\centering
		\begin{tabular}{ c  c  c  c  c  c }
			\hline
			\texttt{lc(p1)} & \texttt{lc(p2)} & \texttt{lc(p3)} & - & - & - \\
			1.2437 & 1.2414 & 1.2404 & - & - & - \\
			\hline
			$\left\| \texttt{p1}-\texttt{x\_opt} \right\| $ & $\left\| \texttt{p2}-\texttt{x\_opt} \right\| $ & $\left\| \texttt{p3}-\texttt{x\_opt} \right\| $ & $\left\| \texttt{lc(p1)}-\texttt{fval} \right\| $ & $\left\| \texttt{lc(p2)}-\texttt{fval} \right\| $ & $\left\| \texttt{lc(p3)}-\texttt{fval} \right\| $\\
			 0.1286 & 0.0835 &  0.0544 & 0.0039 & 0.0015 & 5.6546e-004 \\
			\hline
			\hline
			\texttt{lc(q1)} & \texttt{lc(q2)} & - & - & - & - \\
			1.4686 & 1.5508 & - & - & - & - \\
			\hline
			$\left\| \texttt{p}-\texttt{q1} \right\| $ & $\left\| \texttt{p}-\texttt{q2} \right\| $ & $\left\| \texttt{lc(p)}-\texttt{lc(q1)} \right\| $ & $\left\| \texttt{lc(p)}-\texttt{lc(q2)} \right\| $ & - & -  \\
			0.1414 & 0.1414 & 0.0404 & 0.0418 & - & - \\
			\hline
		\end{tabular}
		\normalsize
		\caption{\small Condition lengths of the B\'{e}zier curves are defined by \texttt{p1, p2, p3}, \texttt{q1, q2}, and the perturbations with their optimum.}\label{tab3.2}
		\normalsize
	\end{table}
	Based on the numerical findings in Table \ref{tab3.2}, we can conclude the optimal value of the condition length does not change much, with the gap values less than $0.0016$. \\
	On the other hand, we are going to perturb the initial control points a bit to see how large the condition length variation is. We consider other sets of control points $\left\lbrace \binom{-1}{-1}, \binom{-0.5}{-1.9}, \binom{0.5}{-2.4},\binom{1}{-1}\right\rbrace $ and $\left\lbrace \binom{-1}{-1},\binom{-0.5}{-2.1},\binom{0.5}{-2.6},\binom{1}{-1}\right\rbrace $, which we denote using matrices \texttt{q1} and \texttt{q2} as below
	\small
	$$\texttt{q1} = \begin{pmatrix}
	-1.0000 & -0.5000  &  0.5000 & 1.0000 \\
	-1.0000 & -1.9000 & -2.4000 & -1.0000
	\end{pmatrix},\quad \texttt{q2} = \begin{pmatrix}
	-1.0000 & -0.5000  &  0.5000 & 1.0000 \\
	-1.0000 & -2.1000 & -2.6000 & -1.0000
	\end{pmatrix}.$$
	\normalsize
	\begin{figure}[!ht]
		\subfloat[Starting control points \texttt{p}, \texttt{q1} and \texttt{q2}, and its B\'{e}zier curve of degree $3$.\label{subfig-1:6}]{%
			\includegraphics[width=0.45\textwidth]{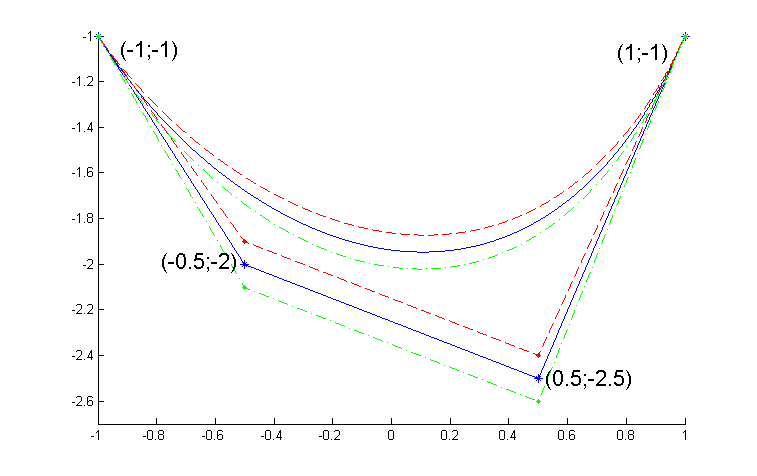}
		}
		\hfill
		\subfloat[Optimal control points \texttt{x\_opt} and its approximate geodesic.\label{subfig-2:6}]{%
			\includegraphics[width=0.45\textwidth]{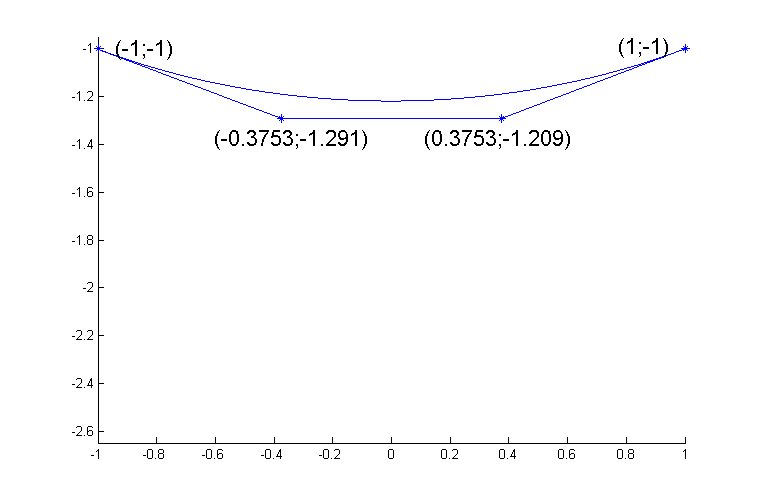}
		}
		\caption{\small B\'{e}zier curve of degree 3 that approximates the condition geodesic of endpoints $\binom{-1}{-1}$ and $\binom{1}{-1}$. Source: Author.}
		\label{fig:6}
	\end{figure}
	\normalsize
	By computing the condition length of the B\'{e}zier curves defined by \texttt{q1} and \texttt{q2} and taking it's perturbation as in table \ref{tab3.2}, we can conclude that the condition length does not change much, with the gap values less than $0.0419$.
\end{example}

Moreover, we jump to higher degrees of B\'{e}zier curve, degree $10$ and degree $20$, in the following example to see how it impacts the values of the minimum condition length. 
\begin{example}\label{exa:4.3}
	Consider a set of control points $\left\lbrace \binom{-1}{-1},\binom{-0.8}{0}, \binom{-0.6}{-2}, \binom{-0.4}{-1},\binom{-0.2}{-3},\binom{0}{0},\binom{0.2}{-1},\binom{0.4}{-2},\binom{0.6}{0},\binom{0.8}{-3},\binom{1}{-1}\right\rbrace $, therein we mean that these control points for the control polygon of a B\'{e}zier curve of degree $10$, which we denote using a matrix:
	\small
	\begin{verbatim}
	p =
	Columns 1 through 9
	-1.0000   -0.8000   -0.6000   -0.4000   -0.2000    0    0.2000    0.4000    0.6000
	-1.0000         0   -2.0000   -1.0000   -3.0000    0   -1.0000   -2.0000         0
	Columns 10 through 11
	0.8000    1.0000
	-3.0000   -1.0000	
	\end{verbatim}
	\normalsize
	Then, after approximate $3$ minutes to compute by Matlab, we obtain the approximate optimal values of the control points \texttt{x\_opt} and the condition length \texttt{fval} of geodesic joining $\binom{-1}{-1}$ and $\binom{1}{-1}$
	\small
	\begin{verbatim}
	fval = 1.0228
	x_opt =
	Columns 1 through 9
	-1.0000   -0.9745    0.1422   -0.7599   -0.1195   -0.1148   -0.4413    0.1192    0.5731
	-1.0000   -1.0126   -1.5427   -0.7061   -1.6850   -0.9115   -1.2465   -1.3693   -1.2046
	Columns 10 through 11
	0.9932    1.0000
	-1.0031   -1.0000
	\end{verbatim}
	\normalsize
	On the other hand, we consider another set of control points \texttt{q} (here we mean that these control points for the control polygon of a B\'{e}zier curve of degree $20$) which we denote using a matrix:
	\small
	\begin{verbatim}
	q =
	Columns 1 through 9
	-1.0000   -0.9000   -0.8000   -0.7000   -0.6000   -0.5000   -0.4000   -0.3000   -0.2000
	-1.0000         0   -2.0000   -3.0000   -4.0000   -2.0000         0   -3.0000   -5.0000
	Columns 10 through 18
	-0.1000         0    0.1000    0.2000    0.3000    0.4000    0.5000    0.6000    0.7000
	-1.0000   -3.0000   -2.0000   -5.0000   -1.0000   -2.0000         0   -2.0000         0
	Columns 19 through 21
	0.8000    0.9000    1.0000
	-3.0000   -2.0000   -1.0000
	\end{verbatim}
	\normalsize
	Then, after approximate $15$ minutes to compute by Matlab, we obtain the approximate optimal values of the control points \texttt{x\_opt} and the condition length \texttt{fval} of geodesic joining $\binom{-1}{-1}$ and $\binom{1}{-1}$
	\small
	\begin{verbatim}
	fval = 0.9764
	x_opt =
	Columns 1 through 9
	-1.0000   -0.6830   -0.8573   -1.2679   -0.2138   -0.2424   -0.6579   -0.5008   -0.1139
	-1.0000   -1.1473   -1.0050   -1.0339   -1.0363   -1.8057   -0.4658   -1.1808   -2.2698
	Columns 10 through 18
	-0.0601   -0.2262   -0.1579    0.2623    0.6582    0.5991    0.2084    0.2630    0.9391
	0.2984   -2.4573   -0.2764   -2.4519    0.2642   -2.4288   -0.2884   -1.7412   -0.8721
	Columns 19 through 21
	0.2065    0.6006    1.0000
	-1.2538   -1.1871   -1.0000	
	\end{verbatim}
	\normalsize
\end{example}
From the experiences in Examples \ref{exa:4.1}, \ref{exa:4.2} and \ref{exa:4.3}, we imply a table with minimum length corresponding to various degrees:
\begin{table}[htb]
	\centering
	\begin{tabular}{c c c c c }
		\hline
		Degree of B\'{e}zier curve & 2 & 3 & 10 & 20 \\
		Minimum condition length & 1.3948 & 1.2398 & 1.0228 & 0.9764 \\
		\hline
	\end{tabular}
	\caption{\small Minimum condition lengths correspond to different degrees in the space of univariate polynomial of degree 2.}\label{tab1}
	\normalsize
\end{table}

\begin{figure}[!ht]
	\subfloat[Control points \& approximations of geodesics.\label{subfig-1:7}]{%
		\includegraphics[width=0.475\textwidth]{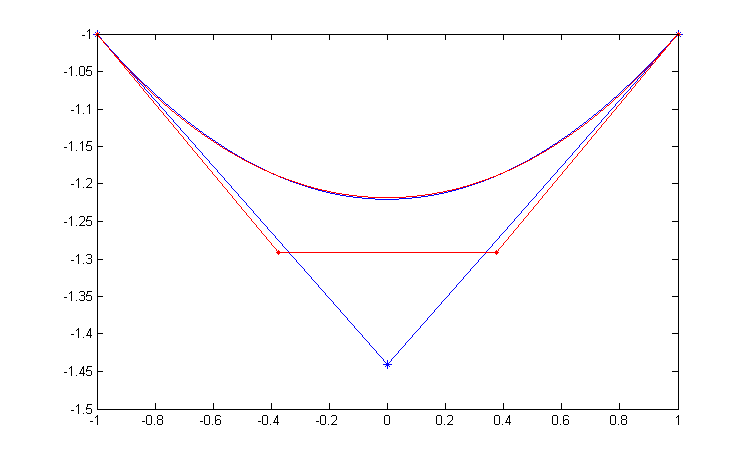}
	}
	\hfill
	\subfloat[Zoom at approximations of  geodesics.\label{subfig-2:7}]{%
		\includegraphics[width=0.475\textwidth]{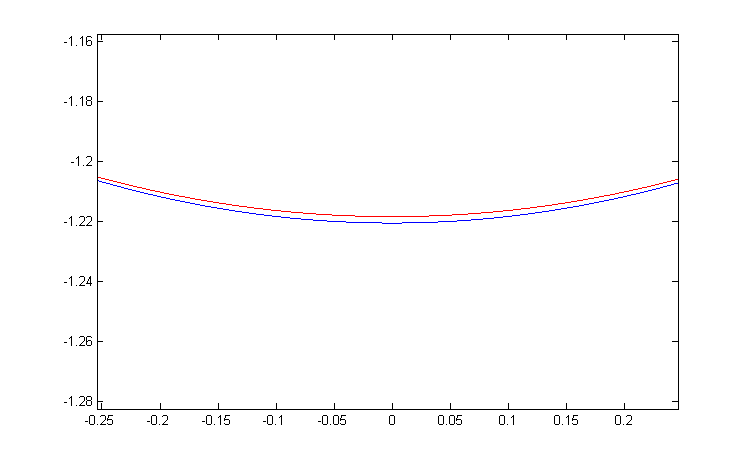}
	}
	\caption{\small Comparison between the approximations by the B\'{e}zier curve degree 2 (in blue) and degree 3 (in red). Source: Author.}
	\label{fig:7}
	\normalsize
\end{figure}
Based on Table \ref{tab1} and Fig. \ref{fig:7}, we can conclude that with more control points, in the space of univariate polynomials of degree $2$, the value of its condition length is better, but above a certain degree the improvement is slight (do not change much).

A more general investigation will take place in the following subsection.
\subsection{In the space of univariate polynomials of degree $n$ ($n\geq 2$)}
In this subsection, \texttt{mlc(p)} indicates the condition length of the B\'{e}zier curve defined by \texttt{p}. We calculate the optimal value, together with the matrix of control points (\texttt{x\_opt}) for the B\'{e}zier curve that realizes the minimum and its condition length (\texttt{fval}).\\
To understand more, we study some examples below, where we use the B\'{e}zier curve defined by these control points as an initial guess for an optimization process that approximates the condition geodesic of endpoints $\left(-1, -1, 2 \right)^T$ and $\left(1, -1, 2 \right)^T$ in the space of univariate polynomials of degree $3$. Note that the endpoints are kept fixed throughout the optimization process.
\begin{example}\label{exa:4.4}
	Consider the set of control points $\left\lbrace \left(-1, -1, 2 \right)^T,\left( 0,-2,2 \right)^T,\left(1, -1, 2 \right)^T\right\rbrace $ which we denote using a matrix:
	\small
	$$ \texttt{p} = \begin{pmatrix}
		-1.0000  &   0   &  1.0000\\
		-1.0000  &  -2.0000   & -1.0000\\
		2.0000   &   2.0000   &  2.0000
	\end{pmatrix}, $$
	\normalsize
	here we consider these control points for the control polygon of a B\'{e}zier curve of degree $2$.\\
	We get the condition length $\texttt{lc(p)} = 0.7622$, and the optimal value
	$$ \texttt{fval} = 0.6341, \quad \texttt{x\_opt} = \begin{pmatrix}
	-1.0000  & -0.0381 & 1.0000 \\
	-1.0000 & -0.9521 & -1.0000 \\
	2.0000 & 2.3306 & 2.0000
	\end{pmatrix}.$$

	Now we are going to perturb the optimal control points a little to see how large is the variation of the optimal condition length. Let's consider other sets of control points\\ $\left\lbrace \left(-1, -1, 2 \right)^T,\left( 0,-0.9521,2.3306 \right)^T,\left(1, -1, 2 \right)^T\right\rbrace $, $\left\lbrace \left(-1, -1, 2 \right)^T,\left( -0.0381,-1,2.3306 \right)^T,\left(1, -1, 2 \right)^T\right\rbrace $ and $\left\lbrace \left(-1, -1, 2 \right)^T,\left( -0.0381,-0.9521,2.2 \right)^T,\left(1, -1, 2 \right)^T\right\rbrace $, which we denote using matrices \texttt{p1, p2} and \texttt{p3} as below:
	\footnotesize
	$$ \texttt{p1} = \begin{pmatrix}
	-1.0000      &   0  &  1.0000 \\         
	-1.0000  & -0.9521  & -1.0000 \\
	2.0000   & 2.3306  &  2.0000 
	\end{pmatrix}, \texttt{p2} = \begin{pmatrix}
	-1.0000      &    -0.0381  &  1.0000 \\         
	-1.0000  & 1.0000  & -1.0000 \\
	2.0000   & 2.3306  &  2.0000 
	\end{pmatrix}, \texttt{p3} = \begin{pmatrix}
	-1.0000      &    -0.0381  &  1.0000 \\         
	-1.0000  & -0.9521  & -1.0000 \\
	2.0000   & 2.2000  &  2.0000 
	\end{pmatrix}. $$
	\normalsize
	\begin{table}[htb]
		\small 
		\centering
		\begin{tabular}{ c c c c}
			\hline
			\texttt{mlc(p1)} & \texttt{mlc(p2)} & \texttt{mlc(p3)} & - \\
			0.6341 & 0.6343 & 0.6359 & - \\
			\hline
			 $\left\| \texttt{p1}-\texttt{x\_opt} \right\| $ & $\left\| \texttt{p2}-\texttt{x\_opt} \right\| $ & $\left\| \texttt{p3}-\texttt{x\_opt} \right\| $ & - \\
			 0.0381 & 0.0479 & 0.1306 & - \\
			 \hline
			  $\left\| \texttt{mlc(p1)}-\texttt{fval} \right\| $ & $\left\| \texttt{mlc(p2)}-\texttt{fval} \right\| $ & $\left\| \texttt{mlc(p3)}-\texttt{fval} \right\| $ & - \\
			9.2920e-006 & 2.4545e-004 & 0.0018 & - \\
			\hline
			\hline
			\texttt{mlc(q1)} & \texttt{mlc(q2)} & \texttt{mlc(q3)} & \texttt{mlc(q4)} \\
			0.7425 & 0.7833 & 0.7605 & 0.7646 \\
			\hline
			 $\left\| \texttt{p}-\texttt{q1} \right\| $ & $\left\| \texttt{p}-\texttt{q2} \right\| $ & $\left\| \texttt{p}-\texttt{q3} \right\| $ & $\left\| \texttt{p}-\texttt{q4} \right\| $ \\
			 0.1000 & 0.1000 & 0.1000 & 0.1000 \\
			\hline
			$\left\| \texttt{mlc(p)}-\texttt{mlc(q1)} \right\| $ & $\left\| \texttt{mlc(p)}-\texttt{mlc(q2)} \right\| $ & $\left\| \texttt{mlc(p)}-\texttt{mlc(q3)} \right\| $ & $\left\| \texttt{mlc(p)}-\texttt{mlc(q4)} \right\| $ \\
			0.0197 & 0.0211	& 0.0017 & 0.0024 \\
			\hline
		\end{tabular}
		\caption{\small Condition lengths of the B\'{e}zier curves are defined by \texttt{p1, p2, p3}, \texttt{q1, q2, q3}, and the perturbations with their optimum.}\label{tab3.4}
		\normalsize
	\end{table}
	By computing the gap values between the condition length of the given control points and the optimal one, see the first half-part of table \ref{tab3.4}, we can conclude that the optimal value of the condition length does not change much, with the gap values less than $0.0019$.\\
	On the other hand, we are going to perturb the initial control points a bit to see how large the condition length variation is. We consider other sets of control points\\
	$\left\lbrace \left(-1, -1, 2 \right)^T,\left( 0,-1.9,2 \right)^T,\left(1, -1, 2 \right)^T\right\rbrace $, $\left\lbrace \left(-1, -1, 2 \right)^T,\left( 0,-2.1,2 \right)^T,\left(1, -1, 2 \right)^T\right\rbrace $,\\ $\left\lbrace \left(-1, -1, 2 \right)^T,\left( 0.1,-2,2 \right)^T,\left(1, -1, 2 \right)^T\right\rbrace $ and $\left\lbrace \left(-1, -1, 2 \right)^T,\left( -0.1,-2,2 \right)^T,\left(1, -1, 2 \right)^T\right\rbrace $, which we denote using matrices \texttt{q1, q2, q3} and \texttt{q4} as follows:
	\small
	$$ \texttt{q1} = \begin{pmatrix}
	-1.0000  &       0  &  1.0000 \\         
	-1.0000  &  -1.9000  & -1.0000 \\
	2.0000   &  2.0000   &  2.0000  
	\end{pmatrix}, \quad \texttt{q2} = \begin{pmatrix}
	-1.0000   &      0    &  1.0000 \\         
	-1.0000   & -2.1000   & -1.0000 \\
	2.0000    & 2.0000    & 2.0000 
	\end{pmatrix}, $$
	$$ \texttt{q3} = \begin{pmatrix}
	-1.0000 &   0.1000 &   1.0000 \\         
	-1.0000 &  -2.0000 &  -1.0000 \\
	2.0000  &  2.0000  &  2.0000  
	\end{pmatrix}, \quad \texttt{q4} = \begin{pmatrix}
	-1.0000 &  -0.1000 &   1.0000 \\
	-1.0000 &  -2.0000 &  -1.0000 \\
	2.0000  &  2.0000  &  2.0000
	\end{pmatrix}. $$
	\normalsize
	The perturbation as in the remain part of table \ref{tab3.4} allows us to conclude that the condition length does not change much, with the gap values less than $0.0212$.
\end{example}
\begin{example}\label{exa:4.5}
	In this example we are going to use a B\'{e}zier curve of degree $3$. Consider the set of control points $\left\lbrace \left(-1,-1,2\right)^T,\left(-0.5,-2,2\right)^T,\left(0.5,-2.5,2\right)^T,\left(1,-1,2\right)^T\right\rbrace $, which is the control polygon of a B\'{e}zier curve of degree $3$, which we denote using a matrix:
	$$ \small \texttt{p} = \begin{pmatrix}
		-1.0000  & -0.5000 &   0.5000 &   1.0000 \\
		-1.0000  & -2.0000 &  -2.5000 &  -1.0000 \\
		2.0000   & 2.0000  &  2.0000  &  2.0000  
	\end{pmatrix},$$
	\normalsize
	The condition length is computed $\texttt{mlc(p)} = 0.8881$. Then we approximate a condition geodesic joining $\left(-1,-1,2 \right)^T$ and $\left(1,-1,2\right)^T$, we obtain the optimal value of the geodesic as follows:
	$$ \texttt{fval} = 0.5635, \quad \texttt{x\_opt} = \begin{pmatrix}
	-1.0000 &  -0.4513 &   0.2635 &   1.0000 \\                  
	-1.0000 &  -0.9453 &  -0.9944 &  -1.0000 \\
	2.0000  &  2.1910  &  2.2461  &  2.0000
	\end{pmatrix}.$$
	Now, to see how large is the variation of the optimal condition length, we are going to perturb the optimal control points a little. We consider other sets of control points \texttt{p1, p2, p3}, and \texttt{p4}, which we denote using matrices as below:
	\small
	$$  \texttt{p1} = \begin{pmatrix}
		-1.0000 &  -0.5513 &   0.1635 &   1.0000 \\
		-1.0000 &  -0.9453 &  -0.9944 &  -1.0000 \\
		2.0000  &  2.1910  &  2.2461  &  2.0000
	\end{pmatrix}, \; \texttt{p2} = \begin{pmatrix}
				-1.0000  & -0.3513  &  0.3635 &   1.0000 \\
				-1.0000  & -0.9453  & -0.9944 &  -1.0000 \\
				2.0000   & 2.1910   & 2.2461  &  2.0000
	\end{pmatrix}, $$
	$$ \texttt{p3} = \begin{pmatrix}
	-1.0000 &  -0.4513 &   0.2635 &   1.0000 \\
	-1.0000 &  -0.8453 &  -0.8944 &  -1.0000 \\
	2.0000  &  2.1910  &  2.2461  &  2.0000
	\end{pmatrix}, \; \texttt{p4} = \begin{pmatrix}
	-1.0000  & -0.4513  &  0.2635 &   1.0000 \\
	-1.0000  & -1.0453  & -1.0044 &  -1.0000 \\
	2.0000   & 2.1910   & 2.2461  &  2.0000
	\end{pmatrix}. $$
	\normalsize
	\begin{table}[htb]
		\small 
		\centering
		\begin{tabular}{c c c c}
			\hline
			\texttt{mlc(p1)} & \texttt{mlc(p2)} & \texttt{mlc(p3)} & \texttt{mlc(p4)} \\
			0.5637 & 0.5636 & 0.5657 & 0.5646 \\
			\hline
			$\left\| \texttt{p1}-\texttt{x\_opt} \right\| $ & $\left\| \texttt{p2}-\texttt{x\_opt} \right\| $ & $\left\| \texttt{p3}-\texttt{x\_opt} \right\| $ & $\left\| \texttt{p4}-\texttt{x\_opt} \right\| $ \\
			0.1414 & 0.1414 & 0.1414 & 0.1005 \\
			\hline
			$\left\| \texttt{mlc(p1)}-\texttt{fval} \right\| $ & $\left\| \texttt{mlc(p2)}-\texttt{fval} \right\| $ & $\left\| \texttt{mlc(p3)}-\texttt{fval} \right\| $ & $\left\| \texttt{mlc(p4)}-\texttt{fval} \right\|$ \\
			1.5350e-004 & 1.3685e-004 & 0.0022 & 0.0011 \\
			\hline
			\hline
			$\left\| \texttt{p}-\texttt{q1} \right\| $ & $\left\| \texttt{p}-\texttt{q2} \right\| $ & $\left\| \texttt{p}-\texttt{q3} \right\| $ & $\left\| \texttt{p}-\texttt{q4} \right\| $ \\
			0.1414 & 0.1414 & 0.1414 & 0.1414\\
			\hline
			$\left\| \texttt{mlc(p)}-\texttt{mlc(q1)} \right\| $ & $\left\| \texttt{mlc(p)}-\texttt{mlc(q2)} \right\| $ & $\left\| \texttt{mlc(p)}-\texttt{mlc(q3)} \right\| $ & $\left\| \texttt{mlc(p)}-\texttt{mlc(q4)} \right\| $ \\
			0.0093 & 0.0431 & 0.0068 & 0.0408 \\
			\hline
		\end{tabular}
		\caption{\small Condition lengths of the B\'{e}zier curves are defined by \texttt{p1, p2, p3}, \texttt{p4}, and the perturbations with their optimum.}\label{tab3.5}
		\normalsize
	\end{table}
	By computing the condition length of the B\'{e}zier curves defined by \texttt{p1, p2, p3}, and \texttt{p4}; and the gap values between the optimal condition length and the condition length of the given control points, see the first half-part in Table \ref{tab3.5}, we can conclude that the optimal value of the condition length does not change much, with the gap values less than $0.0023 $.\\
	On the other hand, we are going to perturb the initial control points a little to see how large is the variation of the condition length. We consider other sets of control points \texttt{q1, q2, q3}, and \texttt{q4}, which we denote using matrices as follows:
	\small
	$$\texttt{q1} = \begin{pmatrix}
	-1.0000 &  -0.6000 &   0.4000 &   1.0000 \\
	-1.0000 &  -2.0000 &  -2.5000 &  -1.0000 \\
	2.0000  &  2.0000  &  2.0000  &  2.0000
	\end{pmatrix}, \; \texttt{q2} = \begin{pmatrix}
	-1.0000  & -0.5000  &  0.5000 &   1.0000 \\
	-1.0000  & -2.1000  & -2.6000 &  -1.0000 \\
	2.0000   & 2.0000   & 2.0000  &  2.0000
	\end{pmatrix}, $$
	$$ \texttt{q3} = \begin{pmatrix}
	-1.0000 &  -0.4000  &  0.6000 &   1.0000 \\
	-1.0000 &   -2.0000 &  -2.5000 &  -1.0000 \\
	2.0000  &  2.0000   & 2.0000  &  2.0000
	\end{pmatrix}, \; \texttt{q4} = \begin{pmatrix}
	-1.0000  & -0.5000  &  0.5000 &   1.0000 \\
	-1.0000  & -1.9000  & -2.4000 &  -1.0000 \\
	2.0000   & 2.0000   & 2.0000  &  2.0000
	\end{pmatrix}. $$
	\normalsize
	The evaluations in the remaining part of Table \ref{tab3.5} allow us to conclude that the condition length does not change much, with the gap values less than $0.0432$.
\end{example}

We furthermore jump to higher degrees of B\'{e}zier curve in the following example to see how it affects the values of the optimal condition length.
\begin{example}\label{exa:4.6}
	In this example, we will use B\'{e}zier curves of higher degrees, $10$ and $20$. Consider a set of control points \texttt{p}, here we mean that these control points for the control polygon of a B\'{e}zier curve of degree $10$, which we denote using a matrix:
	\small
	\begin{verbatim}
	p =
	Columns 1 through 8
	-1.0000    0.6294    0.8116   -0.7460    0.8268    0.2647   -0.8049   -0.4430
	-1.0000    0.9298   -0.6848    0.9412    0.9143   -0.0292    0.6006   -0.7162
	2.0000    2.0000    2.0000    2.0000    2.0000    2.0000    2.0000    2.0000
	Columns 9 through 11
	0.0938    0.9150    1.0000
	-0.1565    0.8315   -1.0000
	2.0000    2.0000    2.0000
	\end{verbatim}
	\normalsize
	Then, after Matlab took approximately $16$ minutes, we obtain the approximate optimal values of the control points \texttt{x\_opt} and the condition length \texttt{fval} of geodesic joining $\left(-1,-1,2\right)^T$ and $\left(1,-1,2\right)^T$
	\small
	\begin{verbatim}
	fval = 0.4649
	x_opt =
	Columns 1 through 8
	-1.0000   -0.2445    0.0509   -1.1190    0.6417    0.3709   -0.3367    0.4317
	-1.0000   -0.9259   -1.0548   -0.8673   -1.0736   -0.9279   -1.0007   -0.9813
	2.0000    2.2647    2.0639    2.1550    2.2810    2.0668    2.1561    2.2539
	Columns 9 through 11
	0.6332    0.5116    1.0000
	-0.9926   -0.9951   -1.0000
	2.0346    2.1632    2.0000
	\end{verbatim}
	\normalsize
	On the other hand, we consider another set of control points \texttt{q} (here we mean that these control points for the control polygon of a B\'{e}zier curve of degree $20$) which denote using a matrix as below:
	\small
	\begin{verbatim}
	q =
	Columns 1 through 8
	-1.0000    0.5844    0.9190    0.3115   -0.9286    0.6983    0.8680    0.3575
	-1.0000   -0.3658    0.9004   -0.9311   -0.1225   -0.2369    0.5310    0.5904
	2.0000    2.0000    2.0000    2.0000    2.0000    2.0000    2.0000    2.0000
	Columns 9 through 16
	0.5155    0.4863   -0.2155    0.3110   -0.6576    0.4121   -0.9363   -0.4462
	-0.6263   -0.0205   -0.1088    0.2926    0.4187    0.5094   -0.4479    0.3594
	2.0000    2.0000    2.0000    2.0000    2.0000    2.0000    2.0000    2.0000
	Columns 17 through 21
	-0.9077   -0.8057    0.6469    0.3897    1.0000
	0.3102   -0.6748   -0.7620   -0.0033   -1.0000
	2.0000    2.0000    2.0000    2.0000    2.0000
	\end{verbatim}
	\normalsize
	Then, after approximate $80$ minutes to compute by Matlab, we get the approximate optimal values of the control points \texttt{x\_opt} and the condition length \texttt{fval} of geodesic joining $\left(-1,-1,2\right)^T$ and $\left(1,-1,2\right)^T$
	\small
	\begin{verbatim}
	fval = 0.4438
	x_opt =
	Columns 1 through 8
	-1.0000   -0.3013   -0.1632   -1.3171    0.6966    0.4633   -0.5298   -0.3239
	-1.0000   -0.9316   -1.0445   -0.8891   -1.0337   -0.9462   -1.0328   -0.9272
	2.0000    2.2428    2.0583    2.0917    2.3671    2.0814    2.0501    2.2592
	Columns 9 through 16
	0.0083    0.7210    0.7007   -0.9750    0.5664    0.5808   -0.0583    1.0567
	-0.9042   -1.2569   -0.5426   -1.5498   -0.3760   -1.4563   -0.7693   -1.0160
	2.2444    2.0851    2.0956    2.2166    2.2046    2.1124    2.1334    2.1553
	Columns 17 through 21
	0.0677    0.4384    0.6902    0.6497    1.0000
	-1.0314   -0.9565   -1.0044   -0.9956   -1.0000
	2.0731    2.2154    2.0572    2.1176    2.0000
	\end{verbatim}
	\normalsize
\end{example}
From the values (\texttt{fval}) in Examples \ref{exa:4.4}, \ref{exa:4.5} and \ref{exa:4.6}, we can take the comparison from the numerical values as in Table \ref{tab2}.

\begin{table}[htb]
	\centering
	\begin{tabular}{c c c c c}
		\hline
		Degree of B\'{e}zier curve & 2 & 3 & 10 & 20 \\
		Minimum condition length & 0.6341 & 0.5635 & 0.4649 & 0.4438 \\
		\hline
	\end{tabular}
	\caption{\small Minimum condition length corresponding to various degrees in the space of univariate polynomial of degree 3.}\label{tab2}
	\normalsize
\end{table}
We conclude that with more control points, in the space of univariate polynomial of degree $3$, the value of its condition length is better, but above a certain degree the improvement is slight, i.e. do not change much.
\begin{remark}\rm
	From Table \ref{tab1} and Table \ref{tab2}, we conclude that when we increase the degree of B\'{e}zier curve and unchanged the dimension of space, the condition length is better (but do not change much), and on the other hand, when we keep stable the degree of B\'{e}zier curve and increase the dimension of space, the improvement of condition length is significant. 
\end{remark}
\subsection{The link between the complexity and the condition length}
We will investigate the link between $k$ (number of steps of the homotopy method) and the condition length. We refer to \cite{D10} for the following presentation.\\
Consider multi-variate polynomial mappings:
\[ f=\left(f_1,\ldots,f_n\right):\mathbb{K}^n\rightarrow \mathbb{K}^n, \]
where $f_i=f_i(z_1,\ldots,z_n)$ is a polynomial with coefficients in $\mathbb{K}$ of degree $d_i$ (here $\mathbb{K}=\mathbb{R}$ or $\mathbb{C}$). We define the associated homogeneous system of $f$ by
\[ F=\left(F_1,\ldots,F_n\right):\mathbb{K}^{n+1}\rightarrow \mathbb{K}^n, \]
with 
\[ F_i\left( z_0,z_1,\ldots,z_n\right)=z_0^{d_i}f\left(\frac{z_1}{z_0},\ldots,\frac{z_n}{z_0} \right). \]
Denote as $\mathcal{H}_d$ the space of the homogeneous systems $F$ with $deg(F_i)=d_i$, $d=\left(d_1,\ldots,d_n\right)$. We denote $D=\max \left\lbrace d_i: 1\leq i \leq n \right\rbrace$.\\
Consider the problems-solutions variety
\[ \mathcal{V}=\left\lbrace (F,z)\in\mathbb{P}\left(\mathcal{H}_d\right)\times \mathbb{P}_n(\mathbb{C})\;:\; F(z)=0 \right\rbrace,  \]
and two associated restriction projections $\Pi_1$ and $\Pi_2$ on the coordinate spaces:
\[ \Pi_1:(F,z)\in \mathbb{P}\left(\mathcal{H}_d\right)\times \mathbb{P}_n\left(\mathbb{C}\right)\rightarrow F\in\mathbb{P}\left(\mathcal{H}_d\right),\]
and
\[ \Pi_2:(F,z)\in \mathbb{P}\left(\mathcal{H}_d\right)\times \mathbb{P}_n\left(\mathbb{C}\right)\rightarrow z\in\mathbb{P}_n\left(\mathbb{C}\right).\]
Let $\Sigma'\subseteq \mathcal{V}$ be the set of all critical points of $\Pi_1$ and $\Sigma:=\Pi_1\left(\Sigma'\right)$ the set of critical values.\\
Given $(F,z)\in\mathcal{V}\setminus \Sigma'$, we see that
\[ D\Pi_1(F,z):T_{(F,z)}\mathcal{V}\rightarrow T_F\mathbb{P}\left(\mathcal{H}_d\right),\]
is an isomorphism, so by the inverse function theorem, we can locally reverse the projection $\Pi_1$. By composition with the projection $\Pi_2$ one obtains the solution application
\[ \mathcal{S}_{(F,z)}=\Pi_2\circ\Pi_1:V_F\subset\mathbb{P}\left(\mathcal{H}_d\right)\rightarrow V_z\subset \mathbb{P}_n\left(\mathbb{C}\right), \]
where $V_F$ and $V_z$ are neighborhoods of $F$ in $\mathbb{P}\left( \mathcal{H}_d\right)$ and of $z$ in $\mathbb{P}_n\left( \mathbb{C}\right)$, respectively.\\
The variations in the first order of $z$ in term the variations of $F$ are described by the derivative $D\mathcal{S}_{(F,z)}(F)$ which is given by
\[ D\mathcal{S}_{(F,z)}(F):T_F\mathbb{P}\left(\mathcal{H}_d\right) \rightarrow T_z\mathbb{P}_n\left( \mathbb{C}\right),\quad \dot{z}=D\mathcal{S}_{(F,z)}(F)(\dot{F})=-\left(DF(z)\mid_{z^{\perp}} \right)^{-1}\dot{F}(z). \]
The condition number of $F$ on $z$ is the norm of the operator $D\mathcal{S}_{(F,z)}(F)$:
\[ \mu\left( F,z\right) =\max\limits_{\dot{F}\in T_F\mathbb{P}\left(\mathcal{H}_d\right)}\frac{\left\|D\mathcal{S}_{(F,z)}(F)(\dot{F})\right\|_z}{\left\| \dot{F}\right\|_F}=\left\| F\right\| \left\|\left(DF(z)\mid_{z^{\perp}}\right)^{-1}\text{Diag}\left(\left\| z\right\|^{d_i-1}\right)\right\|, \]
where the last norm is the operator norm defined on $\mathbb{C}^n$ and $\mathbb{C}^{n+1}$. We also see that $\mu(F,z)=\infty$ when $DF(z)\mid_{z^{\perp}}$.\\
The normalized condition number is a variant of $\mu(F,z)$ defined by
\[ \mu_{norm}(F,z)=\left\| F\right\| \left\| \left( DF(z)\mid_{z^{\perp}}\right)^{-1}\text{Diag}\left(\left\| z\right\|^{d_i-1}d_i^{1/2} \right)\right\| .  \]
Given a system $F_1\in\mathbb{P}\left( \mathcal{H}_d\right)$, we want to find a solution $z_1\in\mathbb{P}_n\left( \mathbb{C}\right)$, the homotopy method consists in including this particular problem into a family $F_t\in\mathbb{P}\left( \mathcal{H}_d\right)$, $0\leq t\leq 1$.
\begin{theorem}\cite[\text{Theorem $3.1$}, \text{p. $9$}]{D10}
	Given a curve $F_t\in\mathbb{P}\left( \mathcal{H}_d\right)\setminus \Sigma$, $0\leq t\leq 1$, of class $C^1$ and a solution $z_0\in\mathbb{P}_n\left( \mathbb{C}\right)$ of $F_0$, there exists a unique curve $z_t\in\mathbb{P}_n\left( \mathbb{C}\right)$ which is $C^1$ and which satisfies $\left(F_t,z_t \right)\in \mathcal{V}$.
\end{theorem}
Within the homogeneous background we have chosen, we see the basic equation
\[ F_t(z_t)=0,\quad F_0 \text{ and } z_0 \text{ given}, \]
is equivalent to the initial condition problem
\[ \frac{d}{dt}F_t(z_t)=\dot{F}_t(z_t)+DF_t(z_t)(\dot{z}_t)=0, \quad z_0 \text{ given}, \]
(here $\dot{F}_t$ and $\dot{z}_t$ are the derivatives with respect to $t$), i.e,
\[ \dot{z}_t=-\left( DF_t(z_t)\mid_{z_t^{\perp}}\right)^{-1} \dot{F}_t(z_t), \quad z_0 \text{ given}. \]
We discretize this equation by replacing the interval $\left[ 0,1\right]$ by the sequence $0=t_0<t_1<\ldots<t_k=1$, the solutions $z_{t_i}$ by the approximations $x_i$ and the derivatives with respect to $t$ by divided differences. One obtains
\[ \frac{x_{i+1}-x_i}{t_{i+1}-t_i}=-\left(DF_{i+1}(x_i)\mid_{x_i^{\perp}} \right)^{-1}\frac{F_{i+1}(x_i)-F_i(x_i)}{t_{i+1}-t_i}, \]
and as $F_i(x_i)$ is close to zero we obtain
\[ x_{i+1}=x_i-\left(DF_{i+1}(x_i)\mid_{x_i^{\perp}} \right)^{-1}F_{i+1}(x_i), \]
which we denote
\[x_{i+1}=N_{F_{i+1}}(x_i).\]
\begin{algorithm}\label{alg:1} The \textbf{prediction-correction algorithm} is stated as follows:
	\begin{itemize}
		\item \textbf{Input}: $F_i$, $0\leq i\leq k$, and $x_0$ with $F_0(x_0)=0$,
		\item \textbf{Iteration}: $x_{i+1}=N_{F_{i+1}}(x_i)$, $1\leq i\leq k-1$,
		\item \textbf{Output}: $x_k$.
	\end{itemize}
\end{algorithm}
The complexity of the \textbf{Algorithm \ref{alg:1}} is measured by the number $k$ of steps necessary to obtain an approximate solution $x_k$ of $F_k$.
\\The following result provides a bound for the number of steps $k$ needed.
\begin{theorem}\label{thm:2}\cite[\text{Theorem $3.4$}, \text{p. $13$}]{D10}
	Given a curve $F_t\in\mathbb{P}(\mathcal{H}_d)\setminus \Sigma$, $0\leq t\leq 1$, a solution $z_0$ of $F_0$ and the corresponding lifted curve $\left( F_t,z_t\right)\in\mathcal{V}\setminus \Sigma'$, there exists a subdivision
	\[ 0=t_0<t_1<\ldots<t_k=1, \]
	such that the sequence $x_i$ built by the above prediction-correction algorithm is made up of approximate zeros of $F_i$ corresponding to solutions $z_i$ and
	\[ k\leq CD^{3/2}\mu_{norm}\left(F,z\right)^2L_F. \]
	$C$ is a universal constant, $L_F$ the length of the curve $F_t$ in $\mathbb{P}\left(\mathcal{H}_d\right)$
	\[ L_F=\int\limits_0^1 \left\| \dot{F}_t\right\|_{F_t}dt, \]
	and
	\[ \mu_{norm}\left(F,z\right)=\sup\limits_{0\leq t\leq 1} \mu_{norm}\left(F_t,z_t\right),\]
	is the condition number of the lifted curve.
\end{theorem}
The condition metric is a natural tool to measure the complexity of a homotopy path method. As
\[ L_{\kappa}\left( F,z\right)\leq \mu_{norm}\left( F,z\right) L\left( F,z\right)\leq \mu_{norm}\left( F,z\right)\mu\left( F,z\right)L_F\leq \mu_{norm}\left( F,z\right)^2L_F, \]
where
\[L(F,z)=\int\limits_0^1 \left\| \frac{d}{dt}\left(F_t,z_t \right) \right\|_{\left(F_t,z_t \right) }dt, \]
with 
\[ \left\|\left(\dot{F},\dot{z} \right) \right\|_{(F,z)}^2=\frac{\left\|\dot{F}\right\|^2}{\left\| F\right\|^2}+\frac{\left\|\dot{z}\right\|^2}{\left\| z\right\|^2}, \quad  L_{\kappa}(F,z)=\int\limits_0^1\left\| \frac{d}{dt}\left(F_t,z_t\right)\right\|_{(F_t,z_t)}\mu_{norm}\left(F_t,z_t\right)dt, \]
we obtain from Theorem \ref{thm:2} a better bound (see general result in \cite[\text{Theorem 3}]{S09}) as follows:
\begin{theorem}\cite[\text{Theorem $4.1$}, \text{p. $16$}]{D10}\label{thm:3}
	Given a curve of class $C^1$, $\left( F_t,z_t\right)\in\mathcal{V}\setminus\Sigma'$, $0\leq t\leq 1$, then
	\[ k\leq CD^{3/2}L_{\kappa}\left( F,z\right) \]
	steps of the prediction-correction algorithm are sufficient to achieve our approximate zero calculation.
\end{theorem}
\begin{remark}\label{rem4}\rm 
	In our application, $D$ is the degree of the polynomial we want to solve and $C$ is a universal constant, so the bound on the number of steps of prediction-correction depends linearly on $L_{\kappa}(F,z)$.\\
	Here we can wonder how the condition length $L_{\kappa}(F,z)$ depends on the degree of the approximation of the geodesic we chose. That is the purpose of the examples below.
\end{remark}
In the following examples, we use the approximation $k\simeq CD^{3/2}l_{cn}(\Gamma)$ from Theorem \ref{thm:3}, where $l_{cn}(\Gamma)$ denotes the condition length of a curve $\Gamma$.
\begin{example}\label{exa:5.7}
	In the space of univariate polynomials of degree $2$, we consider two polynomials $p_1(x)=x^2-x-1$ and $p_2(x)=x^2+x-1$. In term of vectors, $p_1=\binom{-1}{-1}$ and $p_2=\binom{1}{-1}$. We obtain the condition length of the geodesic joining $p_1$ and $p_2$. Combining with the results of Examples \ref{exa:4.1}, \ref{exa:4.2} and \ref{exa:4.3}, we obtain a summary (with $D=2$ and $c_2>0$ is a universal constant) as in Table \ref{tab3}.
	\begin{table}[htb]
		\centering
		\begin{tabular}{ c | c c c }
			\hline
			& B\'{e}zier curve & condition length of geodesic & number of steps\\
			\hline
			& linear & $1.9248$ & $k_1=5.4442c_2$\\
			& quadratic & $1.3948$ & $k_2=3.9451c_2$\\
			$p_1$ & cubic & $1.2398$ & $k_3=3.5067c_2$\\
			\& & degree $4$ & $1.1623$ & $k_4=3.2875c_2$ \\
			$p_2$ & degree $5$ & $1.1158$ & $k_5=3.1560c_2$\\
			& degree $10$ & $1.0228$ & $k_{10}=2.8929c_2$\\
			& degree $20$ & $0.9764$ & $k_{20}=2.7617c_2$\\
			\hline
			\hline
			$p_3$ & linear & $4.9558$ & $k_1=14.0171c_2$\\
			\& & quadratic & $2.4309$ & $k_2=6.8756c_2$\\
			$p_4$ & cubic & $2.15$ & $k_3=6.0811c_2$\\
			\hline
		\end{tabular}
		\caption{\small Condition length of the geodesic joining $p_1$ and $p_2$, $p_3$ and $p_4$, together with the degree of B\'{e}zier curve and the number of steps, respectively.}\label{tab3}
		\normalsize
	\end{table}

	On the other hand, consider two polynomials $p_3(x)=x^2-x-0.1$ and $p_4(x)=x^2+x-0.1$. In term of vectors, $p_3=\binom{-1}{-0.1}$ and $p_4=\binom{1}{-0.1}$. The condition length of the geodesic joining $p_3$ and $p_4$ are calculated in Table \ref{tab3}, with $D=2$ and $c_2>0$ is a universal constant. Observe that in this example, the polynomials are closer to the singular locus by comparing with the previous example.
\end{example}
\begin{remark}\rm
	We have some notes from Table \ref{tab3}.
	\begin{itemize}
		\item[(i)] $\frac{k_{j}}{k_{i}}< 1$, where $j>i$. In other words, higher degree approximations of the geodesic yield a shorter condition length and therefore a smaller number of steps for the homotopy method.
		\item[(ii)] When increasing the degree of B\'{e}zier curve, (the condition length of geodesic is better) we will win a lot the number of steps. 
		\item[(iii)] By comparing the columns "number of steps" in Table \ref{tab3}, we conclude that we need more steps when we take a homotopy path closer to the discriminant.
	\end{itemize}
\end{remark}
\begin{example}\label{exa:5.8}
	In the space of univariate polynomials of degree $3$, consider two polynomials $q_1(x)=x^3-x^2-x+2$ and $q_2(x)=x^3+x^2-x+2$. In term of vectors, $q_1=\left(-1,-1,2 \right)^T $ and $q_2=\left(1,-1,2\right)^T$. Then we obtain the condition length of the geodesic joining $q_1$ and $q_2$. Combine with the results of Examples \ref{exa:4.4}, \ref{exa:4.5} and \ref{exa:4.6}, we imply that $\frac{k_{j}}{k_{i}}< 1$, where $j>i$, and the number of steps decrease significantly  when increasing the degree of the B\'{e}zier curve, see Table \ref{tab5} with $D=3$ and $c_3>0$ is a universal constant.
	\begin{table}
		\centering
		\begin{tabular}{c c c}
			\hline
			B\'{e}zier curve & condition length of geodesic & number of steps\\
			\hline
			linear & $0.8621$ & $k_1=4.4796c_3$\\
			quadratic & $0.6341$ & $k_2=3.2949c_3$\\
			cubic & $0.5635$ & $k_3=2.9280c_3$\\
			degree $4$ & $0.5283$ & $k_4=2.7451c_3$\\
			degree $5$ & $0.5072$ & $k_5=2.6355c_3$\\
			degree $10$ & $0.4649$ & $k_{10}=2.4157c_3$\\
			degree $20$ & $0.4438$ & $k_{20}=2.3060c_3$\\
			\hline
		\end{tabular}
		\caption{\small Condition length of the geodesic joining $q_1$ and $q_2$ together with the degree of B\'{e}zier curve and the number of steps.}\label{tab5}
		\normalsize
	\end{table}
\end{example}
\begin{figure}[!ht]
	\subfloat[Result from Table \ref{tab3} and Example \ref{exa:5.7}.\label{subfig-1:dummy}]{%
		\includegraphics[width=0.45\textwidth]{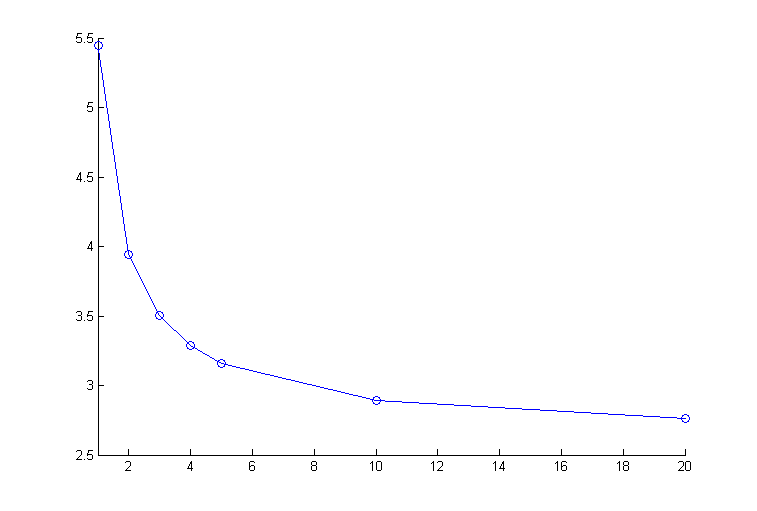}
	}
	\hfill
	\subfloat[Result from Table \ref{tab5} and Example \ref{exa:5.8}.\label{subfig-2:dummy}]{%
		\includegraphics[width=0.45\textwidth]{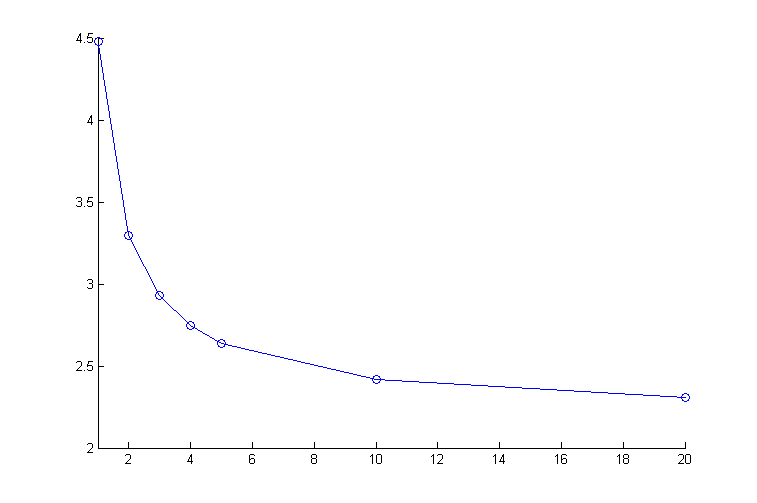}
	}
	\caption{The link between the degree of B\'{e}zier curve and the number of steps. Source: Author.}
	\label{fig:8}
\end{figure}
\begin{remark}\rm 
	Figure \ref{fig:8}, which is obtained from Tables \ref{tab3} and \ref{tab5}, shows the relation between the degree of B\'{e}zier curve and the number of steps: we can conjecture that it is an \textit{exponential} behavior.
\end{remark}
\section{Conclusion}
In this paper, we are concerned with an overview of the path homotopy method for univariate polynomials with the Newton method. The present work defines the condition length of a path joining two polynomials. For these, the B\'{e}zier curves are used to obtain the approximation of geodesics.\\
As a further development, it would be interesting to perform a similar analysis for the case where the correction operator is the Weierstrass method, and all the roots are simultaneously approximating. Moreover, future work may examine the optimization process by the B\'{e}zier surfaces in the space of multivariate polynomials.

\end{document}